\documentclass[a4paper,12pt]{amsart}
\usepackage{amsfonts}
\usepackage{amssymb}
\usepackage{ifthen}
\usepackage{graphicx}
\nonstopmode \numberwithin{equation}{section}
\usepackage[usenames]{color}
\newtheorem{thm}{Theorem}
\newtheorem{lem}{Lemma}
\newtheorem{cor}{Corollary}
\newtheorem{ques}{Question}

\newtheorem{cl}{Claim}
\newtheorem{ca}{Case}
\newtheorem{sca}{Subcase}
\newtheorem{scl}{Subclaim}
\newtheorem{conj}[equation]{Conjecture}

\theoremstyle{definition}
\newtheorem{defn}{Definition}

\newtheorem{op}[equation]{Open Problem}
\newtheorem{rem}{Remark}
\newtheorem{examp}{Example}

\newcounter {own}
\def\theown {\thesection       .\arabic{own}}

\newenvironment{pf}[1][]{%
 \vskip 3mm
 \noindent
 \ifthenelse{\equal{#1}{}}%
  {{\slshape Proof. }}%
  {{\slshape #1.} }%
 }%
{\qed\bigskip}

\newcounter{alphabet}
\newcounter{tmp}
\newenvironment{Thm}[1][]{\refstepcounter{alphabet}%
\bigskip%
\noindent%
{\bf Theorem \Alph{alphabet}}%
\ifthenelse{\equal{#1}{}}{}{ (#1)}%
{\bf .} \itshape}{\vskip 8pt}

\makeatletter
\newcommand{\Ref}[1]{\@ifundefined{r@#1}{}{\setcounter{tmp}{\ref{#1}}\Alph{tmp}}}
\makeatother

\newcommand{\diam}{{\operatorname{Diam}}}

\newcommand{\area}{{\operatorname{Area}}}

\def\be{\begin{equation}}
\def\ee{\end{equation}}

\newcommand{\bee}{\begin{enumerate}}
\newcommand{\eee}{\end{enumerate}}

\newcommand{\blem}{\begin{lem}}
\newcommand{\elem}{\end{lem}}
\newcommand{\bthm}{\begin{thm}}
\newcommand{\ethm}{\end{thm}}
\newcommand{\bcor}{\begin{cor}}
\newcommand{\ecor}{\end{cor}}
\newcommand{\beg}{\begin{exam}}
\newcommand{\eeg}{\end{exam}}
\newcommand{\begs}{\begin{examples}}
\newcommand{\eegs}{\end{examples}}
\newcommand{\bdefe}{\begin{defn}}
\newcommand{\edefe}{\end{defn}}
\newcommand{\bprob}{\begin{prob}}
\newcommand{\eprob}{\end{prob}}
\newcommand{\bques}{\begin{ques}}
\newcommand{\eques}{\end{ques}}
\newcommand{\bei}{\begin{itemize}}
\newcommand{\eei}{\end{itemize}}
\newcommand{\bcon}{\begin{conj}}
\newcommand{\econ}{\end{conj}}
\newcommand{\bop}{\begin{op}}
\newcommand{\eop}{\end{op}}

\newcommand{\bca}{\begin{ca}}
\newcommand{\eca}{\end{ca}}
\newcommand{\bsca}{\begin{sca}}
\newcommand{\esca}{\end{sca}}

\newcommand{\bcl}{\begin{cl}}
\newcommand{\ecl}{\end{cl}}

\newcommand{\bscl}{\begin{scl}}
\newcommand{\escl}{\end{scl}}

\newcommand{\bcons}{\begin{conjs}}
\newcommand{\econs}{\end{conjs}}
\newcommand{\bprop}{\begin{propo}}
\newcommand{\eprop}{\end{propo}}
\newcommand{\br}{\begin{rem}}
\newcommand{\er}{\end{rem}}
\newcommand{\brs}{\begin{rems}}
\newcommand{\ers}{\end{rems}}
\newcommand{\bo}{\begin{obser}}
\newcommand{\eo}{\end{obser}}
\newcommand{\bos}{\begin{obsers}}
\newcommand{\eos}{\end{obsers}}
\newcommand{\bpf}{\begin{pf}}
\newcommand{\epf}{\end{pf}}
\newcommand{\ba}{\begin{array}}
\newcommand{\ea}{\end{array}}
\newcommand{\beq}{\begin{eqnarray}}
\newcommand{\beqq}{\begin{eqnarray*}}
\newcommand{\eeq}{\end{eqnarray}}
\newcommand{\eeqq}{\end{eqnarray*}}

\newcounter{minutes}
\divide\time by 60
\newcounter{hours}
\multiply\time by 60 \addtocounter{minutes}{-\time}
\begin{document}
\bibliographystyle{amsplain}
\title [] {On lengths, areas and Lipschitz continuity of polyharmonic mappings}

\def\thefootnote{}
\footnotetext{ \texttt{\tiny File:~\jobname .tex,
          printed: \number\day-\number\month-\number\year,
          \thehours.\ifnum\theminutes<10{0}\fi\theminutes}
} \makeatletter\def\thefootnote{\@arabic\c@footnote}\makeatother

\author{J. Chen}
\address{J. Chen, Department of Mathematics,
Hunan Normal University, Changsha, Hunan 410081, People's Republic
of China.} \email{jiaolongchen@sina.com}

\author{A. Rasila${}^{~\mathbf{*}}$}
\address{A. Rasila, Department of Mathematics,
Hunan Normal University, Changsha, Hunan 410081, People's Republic
of China, and
Department of Mathematics and Systems Analysis, Aalto University, P. O. Box 11100, FI-00076 Aalto,
 Finland.} \email{antti.rasila@iki.fi}

\author{X. Wang}
\address{X. Wang, Department of Mathematics,
Hunan Normal University, Changsha, Hunan 410081, People's Republic
of China.} \email{xtwang@hunnu.edu.cn}

\keywords{length distortion, area distortion, Lipschitz continuity, polyharmonic mapping \\
${}^{\mathbf{*}}$ Corresponding author}

\begin{abstract}
In this paper, we continue our investigation of polyharmonic mappings in the complex plane. First, we establish two Landau type theorems. We also show a three circles type theorem and an area version of the Schwarz lemma. Finally, we study Lipschitz continuity of polyharmonic mappings with respect to the distance ratio metric.
\end{abstract}

\maketitle \pagestyle{myheadings} \markboth{J. Chen, A. Rasila
and X. Wang}{On lengths, areas and Lipschitz continuity of polyharmonic mappings}

\section{Introduction}\label{csw-sec1}

A complex-valued mapping $F$ in a domain $D$ is
called {\it polyharmonic} (or {\it $p$-harmonic}) if $F$ satisfies the polyharmonic
equation $\Delta^{p}F =\Delta(\Delta^{p-1}F)= 0$ for some $p\in \mathbb{N}^{+}$, where $\Delta$ is the usual complex Laplacian operator
$$\Delta=4\frac{\partial^{2}}{\partial
z\partial\bar{z}}:=\frac{\partial^2}{\partial
x^2}+\frac{\partial^2}{\partial y^2}.$$
In a simply connected domain, a mapping $F$ is polyharmonic if and only if $F$ has
the following representation:
$$F(z) =\sum_{k=1}^{p}|z|^{2(k-1)}G_{k}(z),$$
where each $G_{k}$ is harmonic, i.e., $\Delta G_{k}(z) = 0$ for $k\in \{1,\cdots,p\}$ (cf. \cite{sh1, sh2}). This is also known as the Almansi expansion (see \cite{ar}).
When $p = 1$ (resp. $p=2$), $F$ is harmonic (resp. biharmonic).
The properties of biharmonic mappings have been investigated by many authors (see, e.g., \cite{z1, z2, z3, CW}). We refer to \cite{jt, du}
for the basic properties of harmonic mappings.
References \cite{sh1, sh2, CRW1, CRW2, CRW3, CRW4} are mainly about the geometry properties for certain classes of polyharmonic mappings, such as the starlikeness and convexity, the extremal points, convolution and existence of neighbourhoods. Coefficient estimates, Landau type theorems and radii problems of them are also investigated.

The classical theorem of three circles \cite{a, ro}, also called Hadamard's three circles
theorem, states that if $f$ is an analytic function in the annulus $B(r_{1},r_{2})=\{z:0<r_{1}<|z|=r<r_{2}<\infty\}$,  continuous on $\overline{B(r_{1},r_{2})}$, and $M_{1}$, $M_{2}$ and $M$ are the
maxima of $f$ on the three circles corresponding to $r_{1}$, $r_{2}$ and $r$, respectively, then
$$M^{\log \frac{r_{2}}{r_{1}}}\leq M_{1}^{\log \frac{r_{2}}{r}}M_{2}^{\log \frac{r}{r_{1}}}.$$
If $r_{2}=1$, $M_{1}=r_{1}^{\alpha}$, and $M_{2}=1$, then the Hadamard's result states that, for $r_{1}\leq r=|z|\leq 1$,
$$|f(z)|\leq M_{1}^{\frac{\log r}{\log r_{1}}}=r^{\alpha},$$
where $\alpha$ is an integer.

Let $\mathbb{D}_r$ denote the disk $\{z:|z|<r,z\in \mathbb{C}\}$, and $\mathbb{D}$ the unit disk $\mathbb{D}_1$.
For a polyharmonic mapping $F$, we denote the diameter of the image set of $F(\mathbb{D}_r)$ by $$\diam F(\mathbb{D}_r):=\sup_{z,w\in \mathbb{D}_r}|F(z)-F(w)|$$ (see \cite{b}).

In \cite{poukka}, Poukka proved the following theorem.

\begin{Thm}\label{ThmA}
Suppose $f$ is analytic in $\mathbb{D}$. Then for all positive integers $n$ we have
\be\label{eq4.0}\frac{|f^{(n)}(0)|}{n!}\leq \frac{1}{2}\diam f(\mathbb{D}).\ee
Moreover, equality holds in \eqref{eq4.0} for some $n$ if and only if $f(z)=f(0)+cz^{n}$ for some constant $c$ of modulus $\diam f(\mathbb{D})/2$.
\end{Thm}

We will generalize Theorem A for polyharmonic mappings, and get some coefficient estimates.

For $r\in[0,1)$, the length of the curve $C(r)=\{w=F(re^{i\theta}):\theta\in[0,2\pi]\}$, counting
multiplicity, is defined by
$$l_{F}(r)=\int_{0}^{2\pi}|dF(re^{i\theta})|=r\int_{0}^{2\pi}|F_{z}(re^{i\theta})-e^{-2i\theta}F_{\overline{z}}(re^{i\theta})|d\theta,$$
where $F$ is a polyharmonic mapping defined in $\mathbb{D}$. In particular, let $l_{F}(1)=\sup_{0<r<1}l_{F}(r)$ (cf. \cite{chpr}).
We use the {\it area function} $S_{F}(r)$ of $F$, counting multiplicity, defined by
$$S_{F}(r)=\int_{\mathbb{D}_{r}}J_{F}(z)d\sigma(z),$$
where $d\sigma$ denotes the normalized Lebesgue area measure on $\mathbb{D}$ (cf. \cite{chpr}). In particular,
we let
$$S_{F}(1)=\sup_{0<r<1}S_{F}(r).$$

For a polyharmonic mapping $F$ in $\mathbb{D}$, we use the following standard notations:
$$\lambda_{F}(z)=\min_{0\leq\theta\leq2\pi}|F_{z}(z)
+e^{-2i\theta}F_{\overline{z}}(z)|=\big||F_{z}(z)|-|F_{\overline{z}}(z)| \big|,$$
 $$\Lambda_{F}(z)=\max_{0\leq\theta\leq2\pi}|F_{z}(z)
 +e^{-2i\theta}F_{\overline{z}}(z)|=\big||F_{z}(z)|+|F_{\overline{z}}(z)| \big|.$$
$F$ is said to be {\it $K$-quasiregular}, $K\in[1,\infty)$, if
for $z\in \mathbb{D}$, $\Lambda_{F}(z)\leq K\lambda_{F}(z)$.

In \cite{b}, the authors proved the following area versio of the Schwarz lemma:

\begin{Thm}\label{ThmB}
Suppose $f$ is analytic on the unit disk $\mathbb{D}$. Then the function $\phi_{\area}(r):=(\pi r^{2})^{-1}\area f(\mathbb{D}_{r})$ is strictly increasing for $0<r<1$, except when $f$ is linear, in which case $\phi_{\area}$ is a constant.
\end{Thm}

It is natural to study if a similar area version of the Schwarz lemma exists for polyharmonic mappings.

For a subdomain $G\subset \mathbb{C}$ and for all $z$, $w\in G$, the distance ratio metric $j_{G}$ is defined as
$$j_{G}(z,w)=\log\left(1+\frac{|z-w|}{\min\{d(z,\partial G),d(w,\partial G)\}}\right),$$
where $d(z,\partial G)$ denotes the Euclidean distance from $z$ to $\partial G$. The distance ratio metric was introduced by F. W. Gehring and B. P. Palka \cite{go2} and in the above simplified form by M. Vuorinen \cite{vo2}.
However, neither the quasihyperbolic metric $k_{G}$ nor the distance ratio metric $j_{G}$ are invariant under M\"{o}bius transformations. Therefore, it is natural to ask what the Lipschitz constants are for these metrics under conformal mappings or M\"{o}bius transformations. F. W. Gehring, B. P. Palka and B. G. Osgood proved that these metrics are not changed by more than a factor 2 under M\"{o}bius transformations, see \cite{go1, go2}:

\begin{Thm}\label{ThmC}
 If $G$ and $G'$ are proper subdomain of $\mathbb{R}^{n}$ and if $f$ is a M\"{o}bius transformation of $G$ onto $G'$, then for all $x$, $y\in G$
$$m_{G'}(f(x),f(y))\leq 2m_{G}(x,y),$$
where $m\in \{j,k\}$.
\end{Thm}

In \cite{si, sivw}, the authors considered Lipschitz continuity of the distance-ratio metric under some M\"{o}bius automorphisms of the unit ball and conformal mappings from $\mathbb{D}$ to $\mathbb{D}$.
We will investigate Lipschitz continuity of the distance-ratio metric under certain classes of harmonic mappings.

The organization of this paper is as follows:
In Section \ref{csw-sec3},
we get some coefficient estimates for polyharmonic mappings, and then we obtain two Landau type theorems.
The corresponding results are Theorems \ref{thm3}, \ref{thm4}, \ref{thm8}, and \ref{thm5}. Theorem \ref{thm3} is a generalization of Theorem A, and Theorem \ref{thm4} is a generalization of  \cite[Theorem 3]{chpr}.
In Section \ref{csw-sec2}, we establish a three circles type theorem and an area version of the Schwarz lemma for polyharmonic mappings. Our results are Theorems \ref{thm1} and \ref{thm2}.
Theorem \ref{thm1} is a generalization of \cite[Theorem 1]{chpr}, and Theorem \ref{thm2} is a generalization of Theorem B.
Finally, in Section \ref{csw-sec4}, we investigate the Lipschitz continuity regarding the distance ratio metric for polyharmonic mappings.
The results are Theorems \ref{thm6} and \ref{thm7}. Theorem \ref{thm6} is a generalization of \cite[Theorem 1]{si}.

\section{length of polyharmonic mappings}\label{csw-sec3}

We begin this section by showing useful coefficient estimates for polyharmonic mappings.

\bthm\label{thm3}
Suppose that $F$ is a polyharmonic mapping in $\mathbb{D}$ of the form
\be\label{eq2.1} F(z)=\sum_{n=1}^{p}|z|^{2(n-1)}
\big(h_{n}(z)+\overline{g_{n}(z)}\big)=\sum_{n=1}^{p}|z|^{2(n-1)}
\sum_{j=1}^{\infty}(a_{n,j}z^{j}+\overline{b_{n,j}}\overline{z}^{j}),\ee
 and all its non-zero coefficients $a_{n_{1},j}$, $a_{n_{2},j}$ and
$b_{n_{1},j}$, $b_{n_{2},j}$ satisfy the condition:
\be\label{eq2.8}
\left|\arg \left\{\frac{a_{n_{1},j}}{a_{n_{2},j}} \right\}\right|\leq \frac{\pi}{2}, \;\left|\arg \left\{\frac{b_{n_{1},j}}{b_{n_{2},j}} \right\}\right|\leq\frac{\pi}{2}.\ee
Then
\be\label{eq2.08}\sum_{n=1}^{p}|a_{n,j}|,\;\sum_{n=1}^{p}|b_{n,j}|\leq \frac{\sqrt{p}}{2}\diam F(\mathbb{D}),\ee and
 $$\sum_{n=1}^{p}\left(|a_{n,j}|+|b_{n,j}|\right)\leq \frac{\sqrt{2p}}{2}\diam F(\mathbb{D})$$
for all $n\in\{1,\ldots,p\}$, $j\geq 1$. For $p=1$, the inequalities in \eqref{eq2.08} are sharp for the mappings $F(z)=Cz^{n}$ and $F(z)=C\overline{z}^{n}$, respectively, where $C$ is a constant.

\ethm
\bpf
Let
\begin{align*}
H(z):=&F(z)-F(ze^{i\frac{\pi}{k}})\\
=&\sum_{n=1}^{p}|z|^{2(n-1)}\sum_{j=1}^{\infty}
\left(a_{n,j}z^{j}\big( 1-e^{i\frac{\pi j}{k}}\big) +\overline{b_{n,j}}\overline{z}^{j}\big(1-e^{-i\frac{\pi j}{k}}\big)\right).\\
\end{align*}
Obviously, $|H(z)|\leq \diam F(\mathbb{D})$, and
\begin{align*}
&\frac{1}{2\pi}\int_{0}^{2\pi}|H(z)|^{2}d\theta\\
=&\sum_{1\leq n_{1},n_{2}\leq p} \sum_{j=1}^{\infty}\left(a_{n_{1},j}\overline{a_{n_{2},j}}+b_{n_{1},j}\overline{b_{n_{2},j}}\right)\left|1-e^{i\frac{\pi j}{k}}\right|^{2}r^{2(n_{1}+n_{2}+j-2)}\\
\leq& \diam^{2}F(\mathbb{D}),\\
\end{align*}
where $|z|=r$. Therefore, 
$$\sum_{1\leq n_{1},n_{2}\leq p} \left(a_{n_{1},j}\overline{a_{n_{2},j}}+b_{n_{1},j}\overline{b_{n_{2},j}}\right)\left|1-e^{i\frac{\pi j}{k}}\right|^{2}r^{2(n_{1}+n_{2}+j-2)}
\leq \diam^{2}F(\mathbb{D}),
$$
for all $j\geq 1$. Set $k=j$, and let $r$ tend to 1. Then by the assumption \eqref{eq2.8}, we get
$$\sum_{n=1}^{p} \left(|a_{n,j}|^{2}+|b_{n,j}|^{2}\right)\leq \frac{1}{4}\diam^{2}F(\mathbb{D}).$$
By Cauchy's inequality, we have
$$\sum_{n=1}^{p}|a_{n,j}|,\;\sum_{n=1}^{p}|b_{n,j}|\leq \frac{\sqrt{p}}{2}\diam F(\mathbb{D}),$$ and
$$\sum_{n=1}^{p}\left(|a_{n,j}|+|b_{n,j}|\right)\leq \frac{\sqrt{2p}}{2}\diam F(\mathbb{D})$$
 for all $j\geq 1$. The proof of the theorem is complete.
\epf

\bthm\label{thm4}
Suppose $F$ is a $K$-quasiregular polyharmonic mapping in $\mathbb{D}$ of the form \eqref{eq2.1}, $l_{F}(1) <\infty$, and satisfies the condition:
\be\label{eq2.9}
\left|\arg \left\{\frac{a_{n_{1},j}}{a_{n_{2},j}} \right\}\right|= \left|\arg \left\{\frac{b_{n_{1},j}}{b_{n_{2},j}} \right\}\right|=0,\ee
for non-zero coefficients $a_{n_{1},j}$, $b_{n_{1},j}$, $a_{n_{2},j}$, and $b_{n_{2},j}$. Then for all $n\in\{1,\ldots,p\}$, $j\geq 1$, $$|a_{n,j}|+|b_{n,j}|\leq \frac{Kl_{f}(1)}{2\pi(n+j-1)}.$$
\ethm
\bpf
By a simple computation, we have
$$F_{z}(z)=\sum_{n=1}^{p}\sum_{j=1}^{\infty}
\Big((n+j-1) a_{n,j}z^{n+j-2}\overline{z}^{n-1} +(n-1)\overline{b_{n,j}}z^{n-2}\overline{z}^{n+j-1}\Big),$$
$$F_{\overline{z}}(z)=\sum_{n=1}^{p}\sum_{j=1}^{\infty}
\Big((n-1) a_{n,j}z^{n+j -1}\overline{z}^{n-2} +(n+j-1)\overline{ b_{n,j}}z^{n-1}\overline{z}^{n+j-2} \Big).$$
Then for $j_{0}\geq 1$, we get that
\begin{align*}
&\frac{1}{2\pi}\int_{0}^{2\pi} \frac{F_{z}(z)}{z^{j_{0}-1}} d\theta\\
=&\frac{1}{2\pi}\int_{0}^{2\pi}\sum_{n=1}^{p}\sum_{j=1}^{\infty}
\Big((n+j-1)a_{n,j}z^{n+j-j_{0}-1}\overline{z}^{n-1}  +(n-1)\overline{b_{n,j}}z^{n-j_{0}-1}\overline{z}^{n+j-1}\Big)d\theta\\
=&\sum_{n=1}^{p}(n+j_{0}-1)a_{n,j_{0}}r^{2(n-1)},\\
\end{align*}
and
\begin{align*}
&\frac{1}{2\pi}\int_{0}^{2\pi} \frac{\overline{F_{\overline{z}}(z)}}{z^{j_{0}-1}} d\theta\\
=&\frac{1}{2\pi}\int_{0}^{2\pi}\sum_{n=1}^{p}\sum_{j=1}^{\infty}
\Big((n-1)\overline{a_{n,j}}\overline{z}^{n+j -1}z^{n-j_{0}-1}  +(n+j-1)b_{n,j} z^{n+j-j_{0}-1}\overline{z}^{n-1} \Big)d\theta\\
=&\sum_{n=1}^{p}(n+j_{0}-1)r^{2(n-1)}b_{n,j_{0}},\\
\end{align*}
which give us
\begin{align}\begin{split}\label{eq2.10}
&\left|\sum_{n=1}^{p}(n+j_{0}-1)a_{n,j_{0}}r^{2(n-1)}\right|+\left|\sum_{n=1}^{p}(n+j_{0}-1)b_{n,j_{0}}r^{2(n-1)}\right|\\
=&\left|\frac{1}{2\pi}\int_{0}^{2\pi} \frac{F_{z}(z)}{z^{j_{0}-1}} d\theta\right|+\left|\frac{1}{2\pi}\int_{0}^{2\pi} \frac{\overline{F_{\overline{z}}(z)}}{z^{j_{0}-1}} d\theta\right|\\
\leq &\frac{1}{2\pi}\int_{0}^{2\pi} \frac{\Lambda_{F}(z)}{r^{j_{0}-1}} d\theta.\\
\end{split}\end{align}
It follows from
\begin{align*}
l_{F}(r)=r\int_{0}^{2\pi} |F_{z}(re^{i\theta})-e^{-2i\theta}F_{\overline{z}}(re^{i\theta})| d\theta \geq\frac{r }{K}\int_{0}^{2\pi} \Lambda_{F}(z)d\theta,\\
\end{align*}
that
\be\label{eq2.11}\int_{0}^{2\pi}  \Lambda_{F}(z)d\theta\leq\frac{Kl_{F}(r)}{r}.\ee
\eqref{eq2.10} and \eqref{eq2.11} imply that
$$\left|\sum_{n=1}^{p}(n+j_{0}-1)a_{n,j_{0}}r^{2(n-1)}\right|+\left|\sum_{n=1}^{p}(n+j_{0}-1)b_{n,j_{0}}r^{2(n-1)}\right|\leq \frac{Kl_{F}(r)}{2\pi r^{j_{0}}}.$$
Let $r\rightarrow1^{-1}$. The assumption \eqref{eq2.9} implies
\be\label{eq2.12}\sum_{n=1}^{p}(n+j_{0}-1)(|a_{n,j_{0}}|+|b_{n,j_{0}}|)\leq \frac{Kl_{F}(1)}{2\pi}\ee
for all $j_{0}\geq 1$, and hence,
$$|a_{n,j}|+|b_{n,j}|\leq \frac{Kl_{F}(1)}{2\pi(n+j-1)},$$ for all $k\in \{1,\ldots,p\}$, $j\geq 1$.
The proof of the theorem is complete.

\epf

Next, we establish two Landau type theorems for polyharmonic mappings.

\bthm\label{thm8}
Suppose $F$ is a polyharmonic mapping in $\mathbb{D}$ of the form \eqref{eq2.1}, $\lambda_{F}(0)=\alpha>0$, $\diam F(\mathbb{D})<\infty $, and satisfies the condition
\eqref{eq2.8} for its non-zero coefficients. Then $F$ is univalent in the disk $\mathbb{D}_{r_{0}}$ and $F(\mathbb{D}_{r_{0}})$ contains a univalent disk $\mathbb{D}_{\rho_{0}}$, where $r_{0}$ is the least positive root of the following equation:
$$\alpha- \frac{\sqrt{2p}}{2}\diam F(\mathbb{D})\left(\frac{2r-r^{2}}{(1-r)^{2}}+
\sum_{n=2}^{p}\frac{r^{2(n-1)}}{(1-r)^{2}} +2 \sum_{n=2}^{p}\frac{(n-1)r^{2(n-1)}}{1-r} \right)=0,$$
and
$$\rho_{0}=r_{0}\left(\alpha-\frac{\sqrt{2p}}{2}\diam F(\mathbb{D})\frac{r_{0}}{1-r_{0}}-
\frac{\sqrt{2p}}{2}\diam F(\mathbb{D})\sum_{n=2}^{p}\frac{2r_{0}^{2(n-1)}}{1-r_{0}}\right).$$
\ethm

\bpf The proof of this result is similar to \cite[Theorem 1]{CRW3}, where
$|a_{n,j}|+|b_{n,j}|\leq \frac{\sqrt{2p}}{2}\diam F(\mathbb{D})$ and $\lambda_{F}(0)=\alpha$ is used
instead of $|a_{n,j}|+|b_{n,j}|\leq \sqrt{M^{4}-1}\cdot \lambda_{F}(0)$
for all $(n,j)\not=(1,1)$, and we omit it.
\epf

\begin{figure}
\centering
\includegraphics[width=3.5in]{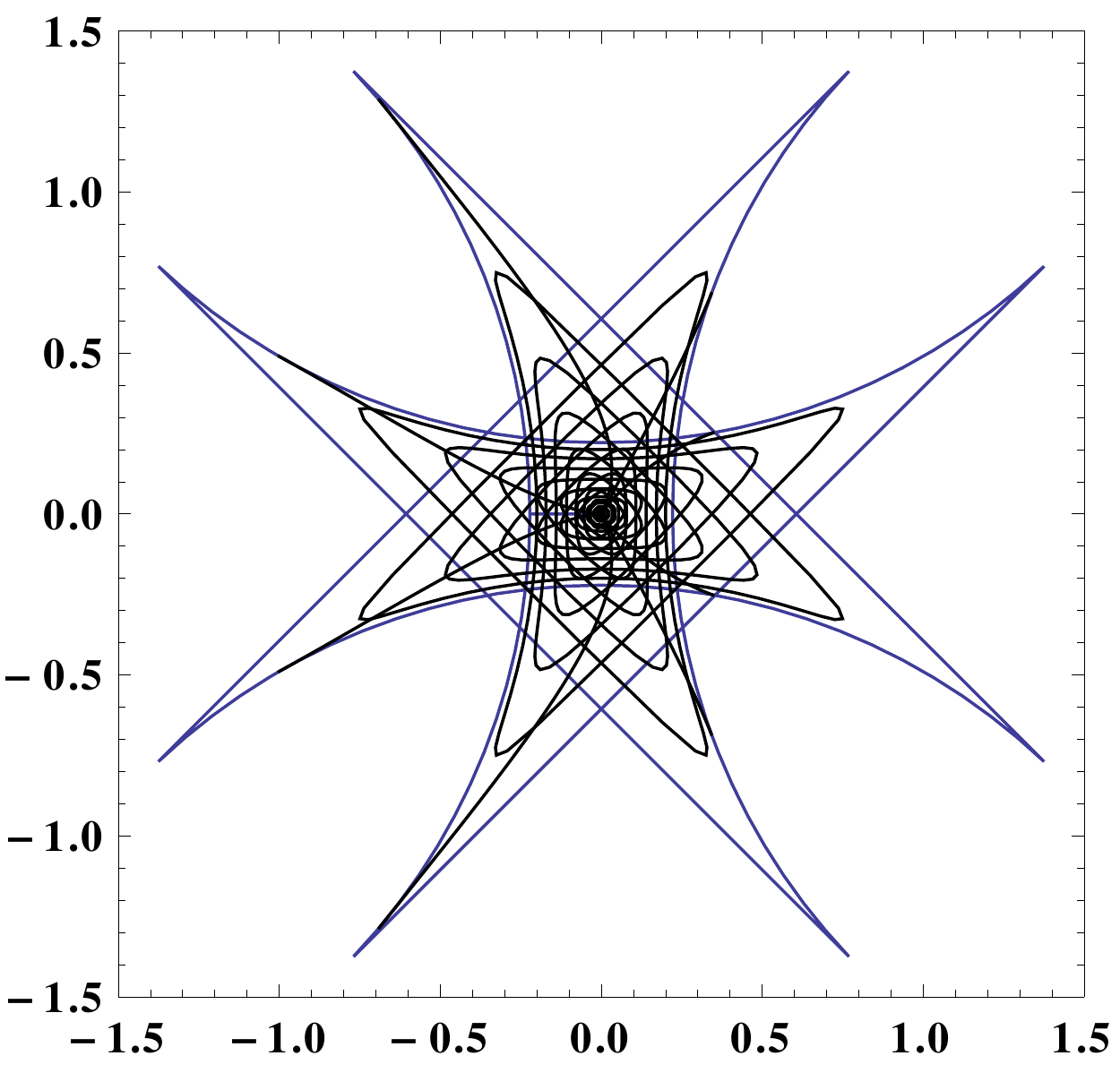}
\caption{The image of the unit disk $\mathbb{D}$ under the mapping $F_{1}$ of Example \ref{ex1}.}\label{F1fig}
\end{figure}

\begin{examp}\label{ex1} Fix $n=4$. Let $\alpha=e^{2\pi i/4}$ be the primitive 4th root of unity, and $\beta=\sqrt{\alpha}=e^{\pi i/4}$. Let $$f_{0}(z)=h_{0}(z)+\overline{ g_{0}(z)}=\frac{1}{\pi}\sum_{k=0}^{3}\alpha^{k}\arg \left\{\frac{z-\beta^{2k+1}}{z-\beta^{2k-1}}\right\}$$
be a harmonic mapping of the disk onto the domain inside a regular 4-gon with vertices at the 4th roots of unity (cf. \cite[p. 59]{du}). By calculations,
$$h_{0}(z)=\sum_{k=0}^{\infty}\frac{4}{\pi (4k+1)}\sin \left(\frac{\pi (4k+1)}{4}\right) z^{4k+1}$$ and $$g_{0}(z)=\sum_{k=1}^{\infty}\frac{4}{\pi (4k-1)}\sin \left(\frac{\pi (4k-1)}{4}\right) z^{4k-1}.$$
Let $F_{1}(z)=\frac{\sqrt{2}\pi}{4}\big(f_{0}(z)+i |z|^{2}f_{0}(z)\big)$ (see Figure \ref{F1fig}). Obviously,  $\lambda_{F_1}(0)=1$, $\diam F_1(\mathbb{D})<\infty$ and the coefficients of $F_{1}$ satisfy the condition
\eqref{eq2.8} for all its non-zero coefficients. Then $F_{1}$ is univalent in the disk $\mathbb{D}_{r_{1}}$ and $F_{1}(\mathbb{D}_{r_{1}})$ contains a univalent disk $\mathbb{D}_{\rho_{1}}$, where $r_{1}$ is the least positive root of the following equation:
$$1-\frac{\sqrt{2p}(r+r^{2}-r^{3})}{(1-r)^{2}}\diam F_1(\mathbb{D})=0,$$
and
$$\rho_{1}=r_{1} \left(
1-\frac{\sqrt{2p} (r_{1}+2r_{1}^{2})}{2(1-r_{1})}\diam F_1 (\mathbb{D})
\right).$$
\end{examp}

\bthm\label{thm5}
Suppose $F$ is a $K$-quasiregular polyharmonic mapping in $\mathbb{D}$ of the form \eqref{eq2.1}, $\lambda_{F}(0)=\alpha>0$, $l_{F}(1) <\infty$, and satisfies the condition
\eqref{eq2.9} for its non-zero coefficients. Then $F$ is univalent in the disk $\mathbb{D}_{r_{2}}$ and $F(\mathbb{D}_{r_{2}})$ contains a univalent disk $\mathbb{D}_{\rho_{2}}$, where $r_{2}$ is the least positive root of equation
$$\alpha-\frac{Kl_{F}(1)}{2\pi(1-r)}\left(r+3\sum_{n=2}^{p}r^{2(n-1)}\right)=0,$$
and $$\rho_{2}=\alpha r_{2}-\frac{Kl_{F}(1)}{2\pi}\left(\log \frac{1}{1-r_{2}}-r_{2}+2\log\frac{1}{1-r_{2}}\sum_{n=2}^{p}r_{2}^{2(n-1)}\right).$$
\ethm

\bpf For any $z_{1}\not=z_{2}$, where $z_{1}$, $z_{2}\in \mathbb{D}_{r}$ and $r\in (0,1)$ is a constant. It follows from \eqref{eq2.12} that
\begin{eqnarray*}
\big |F(z_{1})-F(z_{2})\big |&=&\left |\int_{[z_{1},z_{2}]}F_{z}(z)dz+F_{\overline{z}}(z)d\overline{z}
\right |\\
&\geq &\left |\int_{[z_{1},z_{2}]}F_{z}(0)dz+F_{\overline{z}}(0)d\overline{z}\right |\\
&&- \left |\int_{[z_{1},z_{2}]}\big(F_{z}(z)-F_{z}(0)\big)dz
+\big(F_{\overline{z}}(z)-F_{\overline{z}}(0)\big)d\overline{z}
\right |\\
&\geq& J_1-J_2-J_3-J_4,
\end{eqnarray*}
where
\begin{eqnarray*}
J_1& :=& \left |\int_{[z_{1},z_{2}]}h_{1}'(0)\,dz+
\overline{g_{1}'(0)}d\overline{z}\right |\geq \int_{[z_{1},z_{2}]}\lambda_{F}(0)|dz|
=\lambda_{F}(0)|z_{1}-z_{2}|,\\
J_2&:=& \left |\int_{[z_{1},z_{2}]} \big(h_{1}'(z)-h_{1}'(0)\big )dz
+\big (\overline{g_{1}'(z)}- \overline{g_{1}'(0)}\big )d\overline{z}\right|,\\
&\leq&\int_{[z_{1},z_{2}]}\big|h_{1}'(z)-h_{1}'(0)\big | |dz|
+\big |\overline{g_{1}'(z)}- \overline{g_{1}'(0)}\big| |d\overline{z}|\\
&\leq&|z_{1}-z_{2}|\sum_{j=2}^{\infty}j(|a_{1,j}|+|b_{1,j}|)r^{j-1}\\
&\leq&|z_{1}-z_{2}| \frac{Kl_{F}(1)}{2\pi} \sum_{j=2}^{\infty}r^{j-1}\\
&=&|z_{1}-z_{2}| \frac{Kl_{F}(1)}{2\pi}\frac{r}{1-r},\\\end{eqnarray*}
\begin{eqnarray*}
J_3&  :=& \left |\int_{[z_{1},z_{2}]}\sum_{n=2}^{p}|z|^{2(n-1)}
h'_{n}(z)dz+\sum_{n=2}^{p}|z|^{2(n-1)}\overline{g'_{n}(z)}d\overline{z}\right|,\\
&\leq&\int_{[z_{1},z_{2}]}\sum_{n=2}^{p}|z|^{2(n-1)}
\big|h_{n}'(z)\big| |dz|+\sum_{n=2}^{p}|z|^{2(n-1)} \big|\overline{g_{n}'(z)}\big| |d\overline{z}|\\
&\leq&|z_{1}-z_{2}|\sum_{n=2}^{p}r^{2(n-1)}\sum_{j=1}^{\infty}j(|a_{n,j}|+|b_{n,j}|)r^{j-1}\\
&\leq&|z_{1}-z_{2}|\sum_{n=2}^{p}r^{2(n-1)}\sum_{j=1}^{\infty} \frac{Kl_{F}(1)}{2\pi} r^{j-1}\\
&\leq&|z_{1}-z_{2}|\frac{Kl_{F}(1)}{2\pi}\cdot\sum_{n=2}^{p}\frac{r^{2(n-1)}}{1-r},\\
\end{eqnarray*}
and
\begin{eqnarray*}
J_4& :=& \left |\int_{[z_{1},z_{2}]}\sum_{n=2}^{p}(n-1)|z|^{2(n-2)}
\big(h_{n}(z)+\overline{g_{n}(z)}\big)(\overline{z}dz+zd\overline{z})\right|\\
&\leq&\left |\int_{[z_{1},z_{2}]}2\sum_{n=2}^{p}(n-1)|z|^{2(n-2)}
\sum_{n=1}^{\infty}\left(|a_{n,j}|+|b_{n,j}|\right)|z|^{j+1}|dz|\right|\\
&\leq&2|z_{1}-z_{2}|\sum_{n=2}^{p}r^{2(n-2)}
\sum_{n=1}^{\infty}(n-1)\left(|a_{n,j}|+|b_{n,j}|\right)r^{j+1}\\
&\leq&2|z_{1}-z_{2}|\sum_{n=2}^{p}r^{2(n-2)}
\sum_{n=1}^{\infty}\frac{Kl_{F}(1)}{2\pi}r^{j+1}\\
&=&2|z_{1}-z_{2}|\frac{Kl_{F}(1)}{2\pi}\cdot \sum_{n=2}^{p }\frac{r^{2(n-1)}}{1-r}.
\end{eqnarray*}
That is
$$\big |F(z_{1})-F(z_{2})\big|\geq |z_{1}-z_{2}|\varphi(r),$$
where
$$\varphi(r)=\alpha- \frac{Kl_{F}(1)}{2\pi}\left(\frac{r}{1-r}+3\sum_{n=2}^{p }\frac{r^{2(n-1)}}{1-r}\right),
$$

It is easy to see that the function $\varphi(r)$ is strictly decreasing
for $r\in(0,1)$,
$$\lim_{r\rightarrow 0+}\varphi(r)=\alpha\ \mbox{and}\
\lim_{r\rightarrow 1^{-}}\varphi(r)=-\infty.
$$
Hence there exists a unique $r_{2}\in (0,1)$ satisfying
$\varphi(r_{2})=0.$ This implies that $F $ is univalent in
$\mathbb{D}_{r_{2}}$.

For any $w$ in $\{w:\; |w|=r_{2}\}$, we obtain
\begin{eqnarray*}
\big |F(w)-F(0)\big|
&=&\left |\int_{[0,w]}F_{z}(z)dz+F_{\overline{z}}(z)d\overline{z}
\right |\\
&\geq&\alpha r_{2}-\sum_{j=2}^{\infty}(|a_{1,j}|+|b_{1,j}|)r_{2}^{j}-
2\sum_{n=2}^{p}r_{2}^{2(n-1)}\sum_{j=1}^{\infty}(|a_{n,j}|+|b_{n,j}|)r_{2}^{j}\\
&\geq&\alpha r_{2}- \frac{Kl_{F}(1)}{2\pi}\sum_{j=2}^{\infty}\frac{r_{2}^{j}}{j}-
2\sum_{n=2}^{p}r_{2}^{2(n-1)}\sum_{j=1}^{\infty} \frac{Kl_{F}(1)}{2\pi}\frac{r_{2}^{j}}{j}\\
&=&\alpha r_{2}- \frac{Kl_{F}(1)}{2\pi}\left(\log\frac{1}{1-r_{2}}-r_{2}+2\log\frac{1}{1-r_{2}}\sum_{n=2}^{p}r_{2}^{2(n-1)}\right):=\rho_{2}.\\
\end{eqnarray*}
Obviously,
$$\rho_{2}>r_{2}\left( \alpha- \frac{Kl_{F}(1)}{2\pi}\left(\frac{r_{2}}{1-r_{2}}+3\sum_{n=2}^{p }\frac{r_{2}^{2(n-1)}}{1-r_{2}}\right) \right)=0.$$
The proof of the theorem is complete.
\epf

\begin{examp} Let $F_{2}(z)=z(1+|z|^{2}+|z|^{4})$ be a $K$-quasiregular polyharmonic mapping. Since $$2\left|\frac{\partial F_{2}(z)}{\partial \overline{z}}\right|<\left|\frac{\partial F_{2}(z)}{\partial z}\right|,$$ then we can choose $K=3$.  Obviously,  $\lambda_{F_{2}}(0)=1$, $l_{F_{2}}(1)<\infty$ and the coefficients of $F_{2}$ satisfy the condition
\eqref{eq2.9} for all its non-zero coefficients.  Then $F_{2}$ is univalent in the disk $\mathbb{D}_{r_{3}}$ and $F_{2}(\mathbb{D}_{r_{3}})$ contains a univalent disk $\mathbb{D}_{\rho_{3}}$, where $r_{3}$ is the least positive root of equation
$$1-\frac{Kl_{F_2}(1)}{2\pi(1-r)}\left(r+3r^{2}+3r^{4}\right)=0,$$
and $$\rho_{3}=r_{3}-\frac{Kl_{F_2}(1)}{2\pi}\left(\log \frac{1}{1-r_{3}}-r_{3}+2(r_{3}^{2}+r_{3}^{4})\log\frac{1}{1-r_{3}}\right).$$
\end{examp}

\section{Area distortion under polyharmonic mappings}\label{csw-sec2}

In this section, we investigate the area distortion under polyharmonic mappings. First, we establish a three circles type theorem, involving the area function, for polyharmonic mappings.

\bthm\label{thm1}
Fix $m\in(0,1)$. Suppose that $F$ is a polyharmonic mapping of the form
\eqref{eq2.1},
$S_{F}(r_{1})\leq m$, $S_{F}(1)\leq 1$, $|a_{n,j}|\geq |b_{n,j}|$ for all $n\in\{1,\cdots,p\}$, $j\geq 1$, and all its non-zero coefficients satisfy the condition:
\be\label{eq2.2}
\left|\arg \left\{\frac{a_{n_{1},j}}{a_{n_{2},j}} \right\}\right|\leq \frac{\pi}{2}, \;\left|\arg \left\{\frac{b_{n_{1},j}}{b_{n_{2},j}} \right\}\right|\geq\frac{\pi}{2},\;\text{where}\;n_{1}\not=n_{2}.\ee
Then for $r_{1}\leq r<1$, $S_{F}(r)\leq m^{\frac{\log r}{\log r_{1}}}$.
\ethm

\bpf
By a simple computation, we have
\begin{align*}
\frac{1}{2\pi}\int_{0}^{2\pi}|F_{z}(z)|^{2}d\theta
=&\sum_{1\leq n_{1},n_{2}\leq p}\sum_{j=1}^{\infty}\Big((n_{1}-1)(n_{2}-1)(a_{n_{1},j}\overline{a_{n_{2},j}}+b_{n_{1},j}\overline{b_{n_{2},j}})\\
&+\big(j(n_{1}+n_{2}-2)+j^{2}\big)a_{n_{1},j}\overline{a_{n_{2},j}} \Big)r^{2(n_{1}+n_{2}+j-3)},\\
\end{align*}
and
\begin{align*}
\frac{1}{2\pi}\int_{0}^{2\pi}|F_{\overline{z}}(z)|^{2}d\theta
=&\sum_{1\leq n_{1},n_{2}\leq p}\sum_{j=1}^{\infty}\Big((n_{1}-1)(n_{2}-1)(a_{n_{1},j}\overline{a_{n_{2},j}}+b_{n_{1},j}\overline{b_{n_{2},j}})\\
&+\big(j(n_{1}+n_{2}-2)+j^{2}\big)b_{n_{1},j}\overline{b_{n_{2},j}} \Big)r^{2(n_{1}+n_{2}+j-3)}.\\
\end{align*}
Therefore,
\begin{align}\begin{split}\label{eq2.3}
&\frac{1}{2\pi}\int_{0}^{2\pi}\big(|F_{z}(z)|^{2}-|F_{\overline{z}}(z)|^{2}\big)d\theta\\
=&\sum_{1\leq n_{1},n_{2}\leq p}\sum_{j=1}^{\infty}  j(n_{1}+n_{2}+j-2) \big(a_{n_{1},j}\overline{a_{n_{2},j}}- b_{n_{1},j}\overline{b_{n_{2},j}} \big)r^{2(n_{1}+n_{2}+j-3)}\\
=&\sum_{n=1}^{p}\sum_{j=1}^{\infty} j(2n+j-2)\big(|a_{n,j}|^{2}- |b_{n,j}|^{2} \big)r^{2(2n+j-3)}\\
&+2\sum_{1\leq n_{1}<n_{2}\leq p}\sum_{j=1}^{\infty}  j(n_{1}+n_{2}+j-2)\text{Re}\big(a_{n_{1},j}\overline{a_{n_{2},j}}- b_{n_{1},j}\overline{b_{n_{2},j}} \big)r^{2(n_{1}+n_{2}+j-3)}.\\
\end{split}\end{align}
It follows from the assumption \eqref{eq2.2} that
$$\frac{1}{2\pi}\int_{0}^{2\pi}\big(|F_{z}(z)|^{2}-|F_{\overline{z}}(z)|^{2}\big)d\theta\geq 0,$$
and hence
\begin{align}\begin{split}\label{eq2.4}
S_{F}(r)=&\int_{\mathbb{D}_{r}}J_{F}(z)d\sigma(z)\\
=&\frac{1}{\pi}\int_{0}^{r}\int_{0}^{2\pi}\big(|F_{z}(\rho e^{i\theta})|^{2}-|F_{\overline{z}}(\rho e^{i\theta})|^{2}\big)d\theta \rho d\rho\\
=&\sum_{n=1}^{p}\sum_{j=1}^{\infty}j\big(|a_{n,j}|^{2}- |b_{n,j}|^{2} \big)r^{2(2n+j-2)}\\
&+2\sum_{1\leq n_{1}<n_{2}\leq p}\sum_{j=1}^{\infty}j\text{Re}\big(a_{n_{1},j}\overline{a_{n_{2},j}}- b_{n_{1},j}\overline{b_{n_{2},j}} \big)r^{2(n_{1}+n_{2}+j-2)}\geq0.\\
\end{split}\end{align}
Let $$G(z)=\sum_{1\leq n_{1},n_{2}\leq p}\sum_{j=1}^{\infty}j\big(a_{n_{1},j}\overline{a_{n_{2},j}}- b_{n_{1},j}\overline{b_{n_{2},j}} \big)z^{2(n_{1}+n_{2}+j-2)}.$$ Then the maximum of $G$ on $\mathbb{D}_{r}$ is obtained on the real axis, that is
 $S_{F}(r)=G(r)=\max_{|z|=r}|G(z)|$, where $0<r_{1}\leq r<1$. Hence the result follows from Hadamard's theorem.
As in \cite[Theorem 1]{chpr}, the mapping $F(z)=\alpha z+\beta \overline{z}$, with $|\alpha|^{2}-|\beta|^{2}=1$ shows the sharpness.
\epf

The following theorem is an area version of the Schwarz lemma for polyharmonic mappings.

\bthm\label{thm2}
Suppose that $F$ is a polyharmonic mapping of the form
\eqref{eq2.1}, $|a_{n,j}|\geq |b_{n,j}|$ for all $n\in\{1,\cdots,p\}$, $j\geq 1$, and all its non-zero coefficients satisfy the condition \eqref{eq2.2}.
Then the function $\phi_{\area}(r):=(\pi r^{2})^{-1}\area F(\mathbb{D}_r)$ is strictly increasing for $0<r<1$, except when $F(z)$ has the form \eqref{eq3.7}, in which case $\phi_{\area}$ is a constant.
\ethm

\bpf
It follows from \eqref{eq2.3} that
\begin{align}\begin{split}\label{eq2.5}
&\frac{1}{2\pi}\int_{0}^{2\pi}J_{F}(re^{i\theta})r d\theta\\
=&\sum_{n=1}^{p}\sum_{j=1}^{\infty} j(2n+j-2)\big(|a_{n,j}|^{2}- |b_{n,j}|^{2} \big)r^{2(2n+j-2)-1}\\
&+2\sum_{1\leq n_{1}<n_{2}\leq p}\sum_{j=1}^{\infty}  j(n_{1}+n_{2}+j-2)\text{Re}\big(a_{n_{1},j}\overline{a_{n_{2},j}}- b_{n_{1},j}\overline{b_{n_{2},j}} \big)r^{2(n_{1}+n_{2}+j-2)-1}.\\
\end{split}\end{align}

Let $$A(r):=\area F(\mathbb{D}_r)=\int_{0}^{2\pi}\int_{0}^{r}J_{F}(\rho e^{i\theta})\rho d \rho d\theta.$$
Since $S_{F}(r)=A(r)/\pi$, then the equations \eqref{eq2.4} imply that

\begin{align}\begin{split}\label{eq2.6}
A(r)=&\pi\sum_{n=1}^{p}\sum_{j=1}^{\infty} j\big(|a_{n,j}|^{2}- |b_{n,j}|^{2} \big)r^{2(2n+j-2)}\\
&+2\pi\sum_{1\leq n_{1}<n_{2}\leq p}\sum_{j=1}^{\infty}j\text{Re}\big(a_{n_{1},j}\overline{a_{n_{2},j}}- b_{n_{1},j}\overline{b_{n_{2},j}} \big)r^{2(n_{1}+n_{2}+j-2)}.\\
\end{split}\end{align}
Since
\begin{align*}
\frac{d A(r)}{dr}=&\frac{d}{dr}\int_{0}^{r}\int_{0}^{2\pi}J_{F}(\rho e^{i\theta})\rho d\theta d\rho\\
=& \int_{0}^{2\pi}J_{F}(re^{i\theta})r d\theta,\\
\end{align*}
then by \eqref{eq2.5} and \eqref{eq2.6}, we have
\begin{multline}
\label{eq2.7}
\frac{d A(r)}{dr}-\frac{2A(r)}{r}\\
=2\pi\left(\sum_{n=1}^{p}\sum_{j=1}^{\infty} j(2n+j-3)\big(|a_{n,j}|^{2}- |b_{n,j}|^{2} \big)r^{2(2n+j-2)-1}\right.\\
\left.+2\sum_{1\leq n_{1}<n_{2}\leq p}\sum_{j=1}^{\infty}j(n_{1}+n_{2}+j-3)\text{Re}\big(a_{n_{1},j}\overline{a_{n_{2},j}}- b_{n_{1},j}\overline{b_{n_{2},j}} \big)r^{2(n_{1}+n_{2}+j-2)-1}\right).\\
\end{multline}
By simple calculations and the assumption, we get
\begin{align*}
\frac{d }{dr}\phi_{\area}(r)=&\frac{1}{\pi r^{2}}\frac{d A(r)}{dr}-\frac{2A(r)}{\pi r^{3}} \\
=& \frac{1}{\pi r^{2}}\left( \frac{d A(r)}{dr}-\frac{2A(r)}{r}\right)\geq 0.\\
\end{align*}

If $\phi_{\area}(r)$ is not strictly increasing, then there is $0<s<t<1$,
such that $\phi_{\area}(r)=C$ for every $s\leq r\leq t$. This
implies that $\phi'_{\area}(r)\equiv 0$ on $[s,t]$, then $\frac{d A(r)}{dr}\equiv\frac{2A(r)}{r}$ on $[s,t]$.
By \eqref{eq2.7}, we see $F$ has the following form
 \begin{align}\begin{split}\label{eq3.7}
F(z)=& z\eta e^{i\theta_{1}}+\overline{z}\xi e^{i\varphi_{1}}
+\sum_{k=2}^{\infty}\zeta_{1,k}(z^{k}e^{i\theta_{k}}+\overline{z}^{k}e^{i\varphi_{k}})\\
&+|z|^{2}\sum_{k=1}^{\infty}\zeta_{2,k}\left(z^{k}e^{i\left(\theta_{k}\pm\frac{\pi}{2}\right)}
+\overline{z}^{k}e^{i\left(\varphi_{k}\pm\frac{\pi}{2}\right)}\right),\\
\end{split}\end{align}
where $\eta$, $\xi$, $\zeta_{1,k}$, $\zeta_{2,k}\geq 0$, and $\theta_{k}$, $\varphi_{k}\in \mathbb{R}$.
\epf

\br If $F$ is analytic, then Theorem \ref{thm2} reduces to \cite[Theorem 1.9]{b}, and gives a new proof of it.
\er

Moreover, by \eqref{eq2.7}, we have
$$\lim_{r\rightarrow 0}\phi_{\area}(r)=\lim_{r\rightarrow 0} \frac{\area F(\mathbb{D}_r)}{\pi r^{2}}=J_{F}(0).$$
Therefore the corollary given below follows.
\bcor\label{cor1}
  Suppose that $F$ is a polyharmonic mapping of the form
\eqref{eq2.1}, $|a_{n,j}|\geq |b_{n,j}|$ for all $n\in\{1,\cdots,p\}$, $j\geq 1$,
and all its non-zero coefficients satisfy the condition \eqref{eq2.2}.
If $\area F(\mathbb{D})=\pi$, then
$$\area F(\mathbb{D}_r)\leq \pi r^{2}$$
for every $0<r<1$.
\ecor

\section{Lipschitz continuity of polyharmonic mappings}\label{csw-sec4}

Now, we give a sufficient condition for a polyharmonic mapping to be a contraction, that is to have the Lipschitz constant at most 1.
\bthm\label{thm6}
Let $F(z)$ be a polyharmonic mapping in $\mathbb{D}$ of the form \eqref{eq2.1}. Suppose that there exists a constant $M>0$ such that $F(\mathbb{D})\subset \mathbb{D}_M$ and
\be\label{eq2.13} \sum_{n=1}^{p}\sum_{j=1}^{\infty}(|a_{n,j}|+|b_{n,j}|)\leq M.\ee
Then $$j_{\mathbb{D}_{M}}(F(z),F(w))\leq j_{\mathbb{D}}(z,w).$$
This inequality is sharp.
\ethm
\bpf
For $z,w\in \mathbb{D}$, let's assume that $|F(z)|\geq |F(w)|$ and $0<r=\max\{|z|,|w|\} $.
Since
\begin{align*}
&|F(z)-F(w)|\\
=&\left| \sum_{n=1}^{p}\sum_{j=1}^{\infty}\left(a_{n,j}(|z|^{2(n-1)}z^{j}-|w|^{2(n-1)}w^{j})+\overline{b_{n,j}}(|z|^{2(n-1)}\overline{z}^{j}-|w|^{2(n-1)}\overline{w}^{j})\right)\right|\\
\leq& |z-w|\sum_{n=1}^{p}\sum_{j=1}^{\infty}\frac{\left||z|^{2(n-1)}z^j-|z|^{2(n-1)}w^j+|z|^{2(n-1)}w^j-|w|^{2(n-1)}w^j\right|}{|z-w|}|a_{n,j}| \\
&+ |z-w|\sum_{n=1}^{p}\sum_{j=1}^{\infty}\frac{\left||z|^{2(n-1)}\overline{z}^{j}-|z|^{2(n-1)}\overline{w}^{j}+|z|^{2(n-1)}\overline{w}^{j}-
|w|^{2(n-1)}\overline{w}^{j}\right|}{|z-w|}|b_{n,j}| \\
\leq& |z-w|\sum_{n=1}^{p}\sum_{j=1}^{\infty}\left(|z|^{2(n-1)} \frac{\left|z^{j}-w^{j}\right|}{|z-w|}+|w|^{j} \frac{|z|^{2(n-1)}-|w|^{2(n-1)}}{|z|-|w|}\right)(|a_{n,j}|+|b_{n,j}|)\\
\leq& |z-w|\sum_{n=1}^{p}\sum_{j=1}^{\infty}\left(|z|^{2(n-1)}\sum_{0\leq s+t\leq j-1}|z|^{s}|w|^{t}+|w|^{j}\sum_{0\leq s+t\leq 2n-3}|z|^{s}|w|^{t}\right)(|a_{n,j}|+|b_{n,j}|)\\
\leq &|z-w| \sum_{n=1}^{p}\sum_{j=1}^{\infty}(|a_{n,j}|+|b_{n,j}|)\sum_{s=0}^{2n+j-3}
|z|^{s},\\
\end{align*}
and
\begin{align*}
M-|F(z)|\geq &\sum_{n=1}^{p} \sum_{j=1}^{\infty}\big(|a_{n,j}|+|b_{n,j}|\big)
-\left|\sum_{n=1}^{p}|z|^{2(n-1)}\sum_{j=1}^{\infty}\big(a_{n,j}z^{j}+\overline{b_{n,j}}\overline{z}^{j}\big)\right|\\
\geq& \sum_{n=1}^{p}\sum_{j=1}^{\infty}(|a_{n,j}|+|b_{n,j}|)(1-|z|^{2n+j-2})\\
=&(1-|z|)  \sum_{n=1}^{p}\sum_{j=1}^{\infty}(|a_{n,j}|+|b_{n,j}|)\sum_{i=0}^{2n+j-3}|z|^{i},\\
\end{align*}
then
\begin{align*}
j_{\mathbb{D}_{M}}(F(z),F(w))=&\log\left(1+\frac{|F(z)-F(w)|}{M-|F(z)|}\right)\\
\leq&\log\left(1+\frac{|z-w|  \sum_{n=1}^{p}\sum_{j=1}^{\infty}(|a_{n,j}|+|b_{n,j}|)\sum_{s=0}^{2n+j-3}|z|^{s}}
{(1-|z|)  \sum_{n=1}^{p}\sum_{j=1}^{\infty}(|a_{n,j}|+|b_{n,j}|)\sum_{i=0}^{2n+j-3}|z|^{i}} \right)\\ 
=&\log\left(1+\frac{|z-w|}{1-|z|}\right)\\
\leq&j_{\mathbb{D}}(z,w).
\end{align*}

As the proof in \cite[Theorem 1]{si}, the mapping $F(z)= |z|^{2(p-1)}z^{j}$ or $F(z)= |z|^{2(p-1)}\overline{z}^{j}$ for $p$, $j\geq1$, shows the sharpness.
\epf

In fact, for a harmonic mapping $f(z)$, the condition $|f(z)|<1$ is not sufficient for the inequality \eqref{eq2.13} to hold for the case $M=1$. For example, one may consider the mapping $f(z)=0.26z+0.25\overline{z}+0.25iz^{2}-0.25i\overline{z}^{2}$. Now, we study Lipschitz continuity of harmonic mappings $f$ with respect to the distance ratio metric, without the condition \eqref{eq2.13}.

\bthm\label{thm7}
Let $f(z)= \sum_{j=1}^{p}\big(a_{j}z^{j}+\overline{b_{j}}\overline{z}^{j}\big)$ be a harmonic mapping in $\mathbb{D}$ with $f(\mathbb{D})\subset \mathbb{D}$.
Then $$j_{\mathbb{D}}(f(z),f(w))<\frac{p\sqrt{2p}}{2}\pi j_{\mathbb{D}}(z,w).$$
\ethm
\bpf
Assume that $|f(z)|\geq |f(w)|$ and $r=\max\{|z|,|w|\}$. It follows from Cauchy's inequality and Parseval's relation
$$\sum_{j=1}^{p}(|a_{j}|^{2}+|b_{j}|^{2})=\frac{1}{2\pi}\int_{0}^{2\pi}|f(z)|^{2}\leq 1$$
that
$$\sum_{j=1}^{p}(|a_{j}|+|b_{j}|)\leq \sqrt{2p\sum_{j=1}^{p}(|a_{j}^{2}|+|b_{j}^{2}|)}\leq \sqrt{2p}.$$ Then,
\begin{align*}
|f(z)-f(w)|=&\left|\sum_{j=1}^{p}\left(a_{k}(z^{k}-w^{k})+\overline{b_{k}}(\overline{z}^{k}-\overline{w}^{k})\right)\right|\\
\leq &p|z-w|\sum_{j=1}^{p}\left(|a_{j}|+|b_{j}|\right)\\
\leq&p\sqrt{2p}|z-w|.
\end{align*}
The Schwarz lemma implies that $1-|f(z)|\geq 1-\frac{4}{\pi}\arctan r$.
Therefore,
\begin{align*}
j_{\mathbb{D}}(f(z),f(w))=&\log\left(1+\frac{|f(z)-f(w)|}{1-|f(z)|}\right)\\
\leq&\log\left(1+p\sqrt{2p}\frac{|z-w|}{1-\frac{4}{\pi} \arctan r} \right)\\
=&\log\left(1+p\sqrt{2p}\frac{|z-w|}{1-r}\frac{1-r}{1-\frac{4}{\pi} \arctan r} \right).\\
\end{align*}
Let $\psi(r)=\frac{g(r)}{h(r)}$, where $g(r)=1-r$, $h(r)=1-\frac{4}{\pi} \arctan r$. Since $g(1)=h(1)=0$, $\frac{g'(r)}{h'(r)}=\frac{\pi(1+r^{2})}{4}$ is strictly increasing with respect to $r$, then $\psi(r)$ is increasing from $[0,1)$ onto $[1,\frac{\pi}{2})$.
Hence,
$$j_{\mathbb{D}}(f(z),f(w))<\log\left(1+\frac{p\sqrt{2p}}{2}\pi \frac{|z-w|}{1-r}\right)\leq \frac{p\sqrt{2p}}{2}\pi j_{\mathbb{D}}(z,w).$$
\epf

\section*{Acknowledgements}
\thanks{The research was partly supported by NSF of China and Hunan Provincial Innovation Foundation for Postgraduates (No. 125000-4242), CIMO Scholarship of Finland, CSC of China, and Academy of Finland (No. 269260).
The work was carried out while the first author was visiting Aalto University.}


\begin{thebibliography}{99}

\bibitem{z1} {\sc Z. Abdulhadi} and {\sc Y. Abu Muhanna},
Landau's theorem for biharmonic mappings.
\textit{J. Math. Anal. Appl.} {\bf 338} (2008), 705--709.

\bibitem{z2} {\sc Z. Abdulhadi}, {\sc Y. Abu Muhanna} and {\sc S. Khuri},
On univalent solutions of the biharmonic equation.
\textit{J. Inequal. Appl.} {\bf 5} (2005), 469--478.

\bibitem{z3} {\sc Z. Abdulhadi}, {\sc Y. Abu Muhanna} and {\sc S. Khuri},
On some properties of solutions of the biharmonic equation.
\textit{Appl. Math. Comput.} {\bf 177} (2006), 346--351.

\bibitem{a} {\sc L. V. Ahlfors},
{\it Conformal Invariants, Topics in Geometric Function Theory}.
McGraw-Hill, New York, 1973.


\bibitem{ar} {\sc N. Aronszajn}, {\sc T. Creese} and {\sc L. Lipkin},
{\it Polyharmonic functions}.
Notes taken by Eberhard Gerlach. Oxford Mathematical Monographs. Oxford Science Publications. The Clarendon Press, Oxford University Press, New York, 1983.

\bibitem{b} {\sc R. B. Burckel}, {\sc D. E. Marshall}, {\sc D. Minda}, {\sc P. Poggi-Corradini} and {\sc T. J. Ransford},
Area, capacity and diameter versions of Schwarz's Lemma.
\textit{Conform. Geom. Dyn.} {\bf 12} (2008), 133--152.

\bibitem{chpr} {\sc Sh. Chen}, {\sc S. Ponnusamy} and {\sc A. Rasila},
Lengths, areas and lipschitz-type spaces of planar harmonic mappings.
arXiv:1309.3767v1 [math.CV].


\bibitem{sh1} {\sc Sh. Chen}, {\sc S. Ponnusamy} and {\sc X. Wang},
Bloch and Landau's theorems for planar $p$-harmonic mappings.
\textit{J. Math. Anal. Appl.} {\bf 373} (2011), 102--110.

\bibitem{sh2} {\sc Sh. Chen}, {\sc S. Ponnusamy} and {\sc X. Wang},
Landau's theorem for $p$-harmonic mappings in several complex variables.
\textit{Ann. Polon. Math.} {\bf 103} (2012), 67--87.

\bibitem{CRW1} {\sc J. Chen}, {\sc A. Rasila} and {\sc X. Wang},
On polyharmonic univalent mappings.
\textit{Math. Rep.} {\bf 15} (4) (2013), 343--357.

\bibitem{CRW2} {\sc J. Chen}, {\sc A. Rasila} and {\sc X. Wang},
Starlikeness and convexity of polyharmonic mappings.
\textit{Bull. Belg. Math. Soc. Simon Stevin.} {\bf 21} (2014), 67--82.

\bibitem{CRW3} {\sc J. Chen}, {\sc A. Rasila} and {\sc X. Wang},
 Landau's theorem for polyharmonic mappings.
\textit{J. Math. Anal. Appl.} {\bf 409} (2014), 934--945.

\bibitem{CRW4} {\sc J. Chen}, {\sc A. Rasila} and {\sc X. Wang},
Coefficient estimates and radii problems for certain classes of polyharmonic mappings.
Preprint. 

\bibitem{CW} {\sc J. Chen} and {\sc X. Wang},
On certain classes of biharmonic mappings defined by convolution.
\textit{Abstr. Appl. Anal.}
{\bf 2012}, Article ID 379130, 10 pages. doi:10.1155/2012/379130

\bibitem{jt} {\sc J. G. Clunie} and {\sc T. Sheil-Small},
Harmonic univalent functions.
\textit{Ann. Acad. Sci. Fenn. Ser. A. I.} {\bf 9} (1984), 3--25.

\bibitem{du} {\sc P. Duren},
{\it Harmonic mappings in the plane}.
Cambridge University Press, Cambridge, 2004.

\bibitem{go1} {\sc F. W. Gehring} and {\sc B. G. Osgood},
Uniform domains and the quasihyperbolic metric.
\textit{J. Analyse Math.} {\bf 36} (1979), 50--74.

\bibitem{go2} {\sc F. W. Gehring} and {\sc B.P. Palka},
Quasiconformally homogeneous domains.
\textit{J. Analyse Math.} {\bf 30} (1976), 172--199.


\bibitem{poukka} {\sc K. A. Poukka},
\"{U}ber die gr\"{o}{\ss}te Schwankung einer analytischen Funktion in einem Kreise.
\textit{Arch. der Math. und Physik} {\bf 11} (1907), 302--307.


\bibitem{si} {\sc S. Simi\'{c}},
Lipschitz continuity of the distace ratio metric on the unit disk.
{\it Filomat} {\bf 27:8} (2013), 1505--1509.

\bibitem{sivw} {\sc S. Simi\'{c}}, {\sc M. Vuorinen}, and {\sc G. Wang},
Sharp Lipschitz constants for the distance ratio metric.
In press, \textit{Math. Scand.}

\bibitem{ro} {\sc R. M. Robinson},
Hadamard's three circles theorem.
\textit{Bull. Amer. Math. Soc.} {\bf 50} (1944), 795--802.

\bibitem{vo2} {\sc M. Vuorinen},
Conformal invariants and quasiregular mappings.
\textit{J. Analyse Math.} {\bf 45} (1985), 69--115.


\end{thebibliography}
\end{document}